\RequirePackage{fix-cm}
\documentclass[smallcondensed]{svjour3}     
\smartqed  
\usepackage{graphicx}
\usepackage{cite}
\usepackage[utf8]{inputenc}
\usepackage[T1]{fontenc}
\usepackage{microtype}
\usepackage{caption}
\usepackage{subcaption}
\usepackage{booktabs} 
\usepackage{url}
\usepackage{amsfonts}
\usepackage{nicefrac}
\usepackage{stmaryrd}       
\usepackage{amsmath}
\usepackage{amssymb}
\usepackage{pifont}
\usepackage{dsfont}
\usepackage{enumerate}
\usepackage{xcolor}
\usepackage{pgfplots}
\pgfplotsset{compat=1.15}
\usepackage{mathrsfs}
\usetikzlibrary{arrows}
\usepackage{graphicx}
\usepackage[title]{appendix}
\usepackage{hyperref}

\newcommand{\E}{\mathbb{E}}	
\newcommand{\R}{\mathbb{R}} 
\renewcommand{\P}{\mathbb{P}}
\newcommand{\cmark}{\ding{51}}%
\newcommand{\xmark}{\ding{55}}
\newcommand{\1}{\mathds{1}}
\newcommand{\ind}{\perp\!\!\!\!\perp}

\DeclareUnicodeCharacter{0301}{*************************************}

\definecolor{qqwuqq}{rgb}{0,0.39215686274509803,0}
\definecolor{qqqqff}{rgb}{0,0,1}
\definecolor{ududff}{rgb}{0.30196078431372547,0.30196078431372547,1}

\journalname{Machine Learning}
\begin{document}

\title{Fairness seen as Global Sensitivity Analysis
}


\author{Clément Bénesse        \and
        Fabrice Gamboa         \and
        Jean-Michel Loubes     \and
        Thibaut Boissin 
}


\institute{C. Bénesse \at
              Institut de Mathématiques de Toulouse \\
              \email{clement.benesse@math.univ-toulouse.fr}           
           \and
           F. Gamboa \& J-M. Loubes \at
              Institut de Mathématiques de Toulouse
           \and
           T. Boissin \at
              IRIT Saint-Exupéry, Toulouse
}

\date{Received: date / Accepted: date}

\maketitle

\begin{abstract}
Ensuring that a predictor is not biased against a sensitive feature is the goal of fair learning. Meanwhile, Global Sensitivity Analysis (GSA) is used in numerous contexts to monitor the influence of any feature on an output variable. We merge these two domains, Global Sensitivity Analysis and Fairness, by showing how Fairness can be defined using a special framework based on Global Sensitivity Analysis and how various usual indicators are common between these two fields. We also present new Global Sensitivity Analysis indices, as well as rates of convergence, that are useful as fairness proxies.
\end{abstract}

\section{Introduction}
\label{intro}
Quantifying  the influence of a variable on the outcome of an algorithm is an issue of high importance in order to explain and understand  decisions taken by machine learning models.  In particular, it enables to detect unwanted biases in the decisions that  lead to unfair predictions. This problem has received a growing attention over the last few years in the literature on fair learning for Artificial Intelligence. One of the main difficulty lies in the definition of what is (un)fair and the choices to quantify it.  A large number of measures have been designed to assess algorithmic fairness, detecting whether a model depends on variables, called sensitive variables, that convey an information that is irrelevant for the model, from a legal or a moral point of view. We refer for instance to \cite{dwork2012fairness,  chouldechova2017fair, bookOneto2020} and \cite{del2020review} and references therein for a presentation of different fairness criteria. Most of these definitions stem back to ensuring the independence between a function of an algorithm output and some sensitive feature that may lead to biased treatment. Hence, understanding and measuring the relationships between a sensitive feature $S$, which is typically included in $\mathbf{X}$ or highly correlated to it, and the output of the algorithm $f(\mathbf{X})$ that predicts a target $Y$, enables to detect unfair algorithmic treatments.  Then, ensuring that predictors are  fair is achieved by controlling previous measures, as done  in~\cite{mary2019fairness, williamson2019fairness, HGR_fairness_1, gordaliza2019obtaining, del2020review, chiappa2020general}. If this notion has been extensively studied for classification,  recent work tackle the regression case as in  \cite{HGR_fairness_1,HGR_fairness_2, chzhenPlug-in2020} or \cite{le2020projection}. \\ \indent
Global Sensitivity Analysis (GSA) is used in numerous contexts for quantifying the influence of a set of features on the outcome of a black-box algorithm. Various indicators, usually taking the form of indices between $0$ and $1$, allow the understanding of how much a feature is important. Multiple set of indices have been proposed over the years such as Sobol’ indices, Cramér-von-Mises indices, HSIC -- see \cite{jacques2006sensitivity, daveiga:hal-01128666, iooss2015review, grandjacques2015analyse, gamboa2020global} and references therein. The flexibility in the choice allows for deep understanding in the relationship between a feature and the outcome of an algorithm. While the usual assumption in this field is to suppose the inputs to be independent, some works \cite{jacques2006sensitivity, mara2012variance, grandjacques2015analyse} remove this assumption to go further in the understanding of the possible ways for a feature to be influential. \\ Hence GSA appears to provide a natural framework to understand the impact of sensitive features. This point of view has been considered when using Shapley values  in the context of fairness \cite{fairness_shap}  and thus provide local fairness by explainability. Hereafter we provide a full probabilistic framework to use GSA for fairness quantification in machine learning. \\
\indent Our contribution is two-fold. First, while GSA is usually concerned with independent inputs, we recall extensions of Sobol’ indices to non-independent inputs introduced in \cite{mara2012variance} that offer ways to account for joint contribution and correlations between variables while quantifying the influence of a feature. We propose an extension of Cramér-von-Mises indices based on similar ideas. We also prove the asymptotic normality for these extended Sobol’ indices to estimate them with a confidence interval. Then, we propose a consistent probabilistic framework to apply  GSA's indices  to quantify fairness. We illustrate the strength of this approach by showing that it can model classical fairness criteria, causal-based fairness and new notions such as intersectionality. This provides new conceptual and practical perspectives to fairness in Machine Learning.

The paper is organized as follows. We begin by reviewing existing works on Global Sensitivity Analysis (Section \ref{section:GSA}). We give estimates for the extended Sobol’ and Cramér-von-Mises indices, along with respectively asymptotic normality (Theorem \ref{theo:TCL}). We then present a probabilistic framework for Fairness in which we draw the link between fairness measures and GSA indices, along with applications to causal fairness and intersectional fairness (Section \ref{section:fairness}).

\section{Global Sensitivity Analysis}\label{section:GSA}
The use of complex computer models for the analysis of applications from science or real-life experiments is by now the routine. The models often are expensive to run and it is important to know with as few runs as possible the global influence of one or several inputs on the outcome of the system under study. When the inputs or features are regarded as random elements, and the algorithm or computer code is seen as a black-box, this problem is referred to as Global Sensitivity Analysis (GSA). Note that since we consider the algorithm to be a black-box, we only need the association of an input and its output. This make it easy to derive the influence of a feature for an algorithm for which we do not have access to new runs. We refer the interested reader to \cite{daveiga:hal-01128666} or \cite{iooss2015review} and references therein for a more complete overview of GSA.

The main objective of GSA is to monitor the influence of variables $X_{1},\cdots, X_{p}$ on an output variable, or variable of interest, $f(X)$. For this, we compare, for a feature $X_i$ and the output $f(X)$, the probability distribution $\P_{X_i,f(X)}$ and the product probability distribution $\P_{X_i}\P_{f(X)}$ by using a measure of dissimilarity. If these two probabilities are equal, the feature $X_i$ has no influence on the output of the algorithm. Otherwise, the influence should be quantifiable. For this, we have access to a wide range of indices, generally tailored to be valued in $\left[ 0, 1 \right]$ and sharing a similar property: the greater the index, the greater the influence of the feature over the outcome. Historically, a variance-decomposition -- or Hoeffding decomposition -- is used of the output of the black-box algorithm to have access to a second-order moment metric in the so-called Sobol’ method. However, these methods were originally developed for independent features. For obvious reasons, this framework is not adapted and has limitations in real-life cases. Additionally, Sobol’ methods are intrinsically restrained by the variance-decomposition and others methods have been proposed. We will present two alternatives for Sobol’ indices. The first one solves the issue of non-independent features. The second one circumvents the limitations of working with variance-decomposition. We finish this section by merging these two alternatives, inspired by the works of \cite{azadkia2019simple, gamboa2020global,chatterjee2020new}.

Note that the use of other metrics is common in the GSA literature. Each metric has its own intrinsic advantages and disadvantages which have been extensively studied.  Moreover, independence tests based on these GSA metrics exist, as shown in \cite{meynaoui2019new,gamboa2020global} and techniques such as bootstrap or Monte-Carlo estimates can be used to obtain confidence intervals for such tests. We restrain ourselves to the Sobol’ and Cramér-von-Mises indices because they are historically the basis of GSA literature, computationally tractable and allow for better understanding of usual fairness proxies, as we will show in Section \ref{section:fairness}. We also prove asymptotic normality for extended Sobol’ indices, which is a first to the best of our knowledge.

\subsection{Sobol’ indices}
A popular and useful tool to quantify the influence of a feature on the output of an algorithm are the Sobol’ indices. Initially introduced in \cite{sobol1990sensitivity}, these indices compare, thanks to the Hoeffding decomposition \cite{van2000asymptotic}, the conditional variance of the output knowing some of the input variables with respect to the overall total variance of the output. Such indices have been extensively studied for computer code experiments.

Suppose that we have the relation $ f(\mathbf{X}) = f(X_1, \cdots, X_p)$ where $f$ is a square-integrable algorithm considered as a black-box and $X_1,\cdots,X_p$ inputs, with $p$ the number of features. We denote by $p_\mathbf{X}$ the distribution of $\mathbf{X}$. For now, we suppose the different inputs to be independent, meaning that $p_\mathbf{X} = \otimes_{i=1}^p p_{X_k}$. Then, we can use the Hoeffding decomposition \cite{van2000asymptotic} on $f(\mathbf{X})$ -- sometimes also called ANOVA-decomposition -- so that we may write
\begin{equation}\label{Hoeffding_function}
    f(\mathbf{X}) = \sum_{s \subseteq \llbracket1,p\rrbracket} f_s(X_s),
\end{equation}
where $f_s$ are square-integrable functions and $X_s$ the set $\{X_i, i \in s\}$. We can either assume that $f$ is centered or that $s$ can be the null set in this sum: it does not change anything since we are interested in the variance afterwards. We will consider $V := \mbox{Var}(f(\mathbf{X}))$ and $V_s := \mbox{Var}(f_s(\mathbf{X}_s))$. Note that the elements of the previous sum are orthogonal in the $L^2(p_{\mathbf{X}})$ sense. So, to compute the variance, we can compute it term by term, and obtain
\begin{equation}\label{Hoeffding_variance}
    V = \sum_{k = 1}^p V_k + \sum_{k_2>k_1}^p V_{k_1,k_2} + \cdots + V_{1,\cdots,p}.
\end{equation}
This equation means that the total variance of the output, which is denoted by $V$, can be split into various components that can be readily interpreted. For instance, $V_1$ represents the variance of the output $f(\mathbf{X})$ that is only due to the variable $X_1$ -- that is, how much $f(\mathbf{X})$ will change if we take different values for $X_1$. Similarly, $V_{1,2}$ represents the variance of the output $Y$ that is only due to the combined effect of the variables $X_1$ and $X_2$ once the main effects of each variable has been removed -- that is, how much $f(\mathbf{X})$ will change if we take different values simultaneously for $X_1$ and $X_2$ and remove the changes due to main effects from $X_1$ only or $X_2$ only.

By dividing the $V_{(m)}$ by $V$, with $(m) \subset \llbracket 1,  p\rrbracket$, we obtain:
\begin{equation}
    S_{(m)}:= \frac{V_{(m)}}{V},
\end{equation}which is the expression of the so-called Sobol’ sensitivity indices. The index $S_k$ quantifies the proportion of the output’s variance caused by the input $X_k$ on its own. The index $S_{(m)},k \in (m)$ quantifies the proportion of the output’s variance caused by the input $X_k$ conjointly with other inputs, and is usually called the Total Sobol’ index of $X_k$.

\subsection{Sobol’ indices for non-independent inputs}

In the classic Sobol’ analysis, for an input $f(\mathbf{X})$, two indices , namely the first order and total indices, quantify the influence of the considered feature on the output of the algorithm. When the inputs are not independent, we need to duplicate each index in order to distinguish whether influences caused by correlations between inputs are taken into account or not. Introduced in this framework by \cite{mara2012variance}, we use the Lévy-Rosemblatt theorem to create two mappings of interest.
We denote by $\sim i$ every index other than $i$. We create $2p$ mappings between $p$ independent uniform random variables $U$ and the variables $\mathbf{X}$ either by mapping $ p_{U_1}p_{U_{\sim 1}}$ to $p_{X_i}p_{X_{\sim i}|X_i}$  -- in this case $U_1$ is denoted by $U^i_1$ -- or by mapping $p_{U_{\sim p}} p_{U_{p}}$ to $p_{X_{\sim i}}p_{X_i}$ -- in this case, $U_{\sim p}$ is denoted $U^{i+1}_{\sim p}$.
In the Appendix \ref{Appendix:rosemblatt}, more in-depth details are given.
In the analysis of the influence of an input $X_i$, the first mapping captures the intrinsic influence of other inputs while the second mapping excludes these influences and shows the variations induced by $X_i$ on its own. Each of these two mappings leads to two indices corresponding to classical Sobol’ and Total Sobol’ indices. The influence of every input $X_i$ is therefore represented by four indices, see Table \ref{table:Sobol_indices}.

Hence, the four Sobol’ indices for each variable $X_i , i \in \llbracket1,p \rrbracket$  are defined as followed:
\begin{equation}\label{eq:defSobol’2}
    Sob_i :=  \frac{\mbox{Var}[\E[f(\mathbf{X})|X_i]]}{\mbox{Var}[f(\mathbf{X})]}
\end{equation}
\begin{equation}
\label{eq:defSobol’T2}
    SobT_i  :=  \frac{\E[\mbox{Var}[f(\mathbf{X})|Z_i]]}{\mbox{Var}[f(\mathbf{X})]}
\end{equation}
\begin{equation}\label{eq:defSobol’ind}
    Sob_i^{ind} :=  \frac{\mbox{Var}[\E[f(\mathbf{X})|Z_i]]}{\mbox{Var}[f(\mathbf{X})]}
\end{equation}
\begin{equation}\label{eq:defSobol’Tind}
    SobT_i^{ind}  :=  \frac{\E[\mbox{Var}[f(\mathbf{X})|X_{\sim i}]]}{\mbox{Var}[f(\mathbf{X})]},
\end{equation}

where the random variable $Z_i$ has the distribution $p_{X_i|X_{\sim i}}$ and is equal to $F^{-1}_{X_i|X_{\sim i}}(U^{i+1}_p)$. \vskip .1in

Note that these definitions can be extended to multidimensional variables and thus enabling to consider groups of inputs by replacing the subset $\{i\}$ by a subset $s \subset \{1,\cdots,p\}$ in the formulas.
\begin{remark}
    \label{remark:sobol=sobolind}
  If the features are independent, then for all $i\in \llbracket 1, \cdots, p\rrbracket$, $Sob^{ind}_i=Sob_i$ and $SobT^{ind}_i = SobT_i$. The proof comes from the fact that in the independent case, we have  $U^i_1 = U^{i+1}_p$. 
\end{remark}

\begin{remark}\label{rk:ineq_Sobol}
All previous indices satisfy the following bounds. For all $i \in \{1,\cdots,p\}$, $$0 \leq Sob^{ind}_i \leq Sob_i \leq SobT_i \leq 1 \quad {\rm and} \quad 0 \leq Sob^{ind}_i \leq SobT^{ind}_i \leq SobT_i \leq 1.$$
    We refer to~\cite{mara2012variance} and  to the law of total variance for the proof.  Note that, in general, there are no inequalities between $Sob_i$ and $SobT^{ind}_i$.
\end{remark}

Sobol indices enable to quantify  three typical ways for a feature to modify the output of an algorithm. 
\begin{enumerate}
    \item Direct contribution. Firstly, a variable can be of interest, all by itself, without any correlation or joint contribution with the other variables. Consider for example the case where $f(\mathbf{x}) = x_1 + x_2$ and $x_1$ independent to the rest of the variables. In this example, we would have $Sob_1 = SobT_1 = Sob^{ind}_1= SobT^{ind}_1 = 0.5$, which means that $50\%$ of the variability of the algorithm is caused by the first variable. In this case, the first variable has a non-null impact on its own on the outcome of the algorithm $f$.
    \item Bouncing contribution. A variable can interact with other variables and influence the output only by its impact on the law of the other variables. For example, consider $(x_1,x_2)$ where $x_2 = \alpha x_1 + \varepsilon$ -- where $\varepsilon$ is a centered white noise of variance $\sigma^2$ -- and $f(\mathbf{x}) = x_2$. Then we get $Sob_1 = SobT_1 = (\alpha^2 V(x_1))/(\alpha^2 V(x_1) + \sigma^2)$ while $Sob^{ind}_1 = SobT^{ind}_1 = 0$. The first variable can be highly influent on the outcome of the algorithm $f$, even if it is not directly responsible for these variations. We call this type of interaction a "bouncing effect" since the variable will need to use another input to reach the outcome of the algorithm.
    \item Joint contribution. Lastly, a variable can contribute to an output jointly with other variables. Take for instance the case where $(x_1,x_2)$ are independent and $f(\mathbf{x}) = x_1 \times x_2$. In this case, $Sob_1 = Sob^{ind}_1 = 0 = Sob_2 = Sob^{ind}_2 $ while $SobT_1 = SobT^{ind}_1 = 1 = SobT_2 = SobT^{ind}_2$ This effect is different of the previous one as the distributions of the input variables are independent but their impact is intertwined. In such a case, the effect is visible and measurable by a variation between first-order and total indices.
\end{enumerate}
\begin{table}[t]
\caption{Sobol’ indices: what is taken into account and what is not.}
\vskip 0.15in
\centering
\begin{sc}
\begin{tabular}{ |c|c|c|  }
 \hline
 \multicolumn{3}{|c|}{Sobol’ indices} \\
 \hline
 &  Correlation between & Joint\\
 &  variables & contributions\\
 \hline
 $Sob_i$ & \cmark & \xmark\\
 \hline
 $SobT_i$ & \cmark & \cmark\\
 \hline
 $Sob^{ind}_i$ & \xmark &\xmark\\
 \hline
 $SobT^{ind}_i$ &\xmark & \cmark\\
 \hline
\end{tabular}
\end{sc}
\label{table:Sobol_indices}
\end{table}

These main differences point out why we need four indices in order to assess the sensitivity of a system to a feature. Table \ref{table:Sobol_indices} sums up which index takes correlations or joint contributions into account. The difference between these different indices can be very informative. For example, if the gap between $Sob_i$ and $SobT_i$ or between $Sob^{ind}_i$ and $SobT^{ind}_i$ is big, then the feature $X_i$ is mainly influential because of its joint contributions with the other features on the output. Conversely, if the gap between $Sob^{ind}_i$ and $Sob_i$ or between $SobT^{ind}_i$  and $SobT_i$ is big, a large part of the influence of the feature $X_i$ will be through its intrinsic influence on other features.

These indices can be rewritten as follow, by using the Lévy-Rosemblatt theorem:
\begin{equation}\label{eq:lrdefSobol’2}
    Sob_i := \frac{\mbox{Var}[\E[g_i(\mathbf{U}^i)|U^i_1]]}{\mbox{Var}[g_i(\mathbf{U}^i)]}  
\end{equation}
\begin{equation}
\label{eq:lrdefSobol’T2}
    SobT_i  := \frac{\E[\mbox{Var}[g_{i}(\mathbf{U}^{i})|U^{i}_{\sim 1}]]}{\mbox{Var}[g_{i}(\mathbf{U}^{i})]}
\end{equation}
\begin{equation}\label{eq:lrdefSobol’ind}
    Sob_i^{ind} := \frac{\mbox{Var}[\E[g_{i+1}(\mathbf{U}^{i+1})|U^{i+1}_p]]}{\mbox{Var}[g_{i+1}(\mathbf{U}^{i+1})]}  
\end{equation}
\begin{equation}\label{eq:lrdefSobol’Tind}
    SobT_i^{ind}  := \frac{\E[\mbox{Var}[g_{i+1}(\mathbf{U}^{i+1})|U^{i+1}_{\sim p}]]}{\mbox{Var}[g_{i+1}(\mathbf{U}^{i+1})]} ,
\end{equation}
as explained in detail in \cite{mara2012variance,mara2015non} or in Appendix \ref{Appendix:estim_mara}. 
Monte-Carlo estimation of the extended Sobol’ indices can be computed by using this definitions. These estimators are consistent and converge to the quantities defined as the Sobol’ and independent Sobol’ indices earlier. Additionally, if we write each of these estimates as $A_n/B_n$, we can use the Delta-method theorem to prove a central limit theorem.
\begin{theorem}\label{theo:TCL}
Each index $\mathcal{S}$ in the equations $\eqref{eq:defSobol’2}$ to $\eqref{eq:defSobol’Tind}$ can be estimated  by its empirical counter part $\mathcal{S}_n$ such that:
\begin{enumerate}[(i)]
  \item $\mathcal{S}_n \xrightarrow{a.s} \mathcal{S}$.
  \item $\sqrt{n}(\mathcal{S}_n - \mathcal{S}) \xrightarrow{D} \mathcal{N}(0, \sigma^2)$, with $\sigma^2$ depending on which index we study, see Appendix \ref{Appendix:estim_mara}.
\end{enumerate}
\end{theorem}

\subsection{Cramér-von-Mises indices}
Sobol’ indices are based on a decomposition of the variance, and therefore only quantify influence of the inputs on the second-order moment of the outcome. Many other criteria to compare the conditional distribution of the output knowing some of the inputs to the distribution of the output have been proposed -- by means of divergences, or measures of dissimilarity between distributions for example. We recall here the definition of Cramér-von-Mises indices \cite{gamboa2020global},  an answer to this lack of distributional information that will be of use later in a fairness framework -- see Section \ref{section:fairness}.

\subsubsection{Classical Cramér-von-Mises indices}
The Cramér-von-Mises indices are based on the whole distribution of $f(\mathbf{X})$. They are defined (see \cite{gamboa2020global}), for every input $i$, as follow:
\begin{equation}
   CVM_i := \frac{\int_{\R} \E \left[(\mu(t) - \mu^i(t))^2\right]d\mu(t)}{\int_{\R} \mu(t)(1-\mu(t))d\mu(t)}
\end{equation}
where $\mu(t) := \E\left[\1_{f(\mathbf{X})\leq t}\right]$ is the cumulative distribution function of $Y$ and $\mu^i$ its conditional version $\mu^i(t) := \E\left[\1_{f(\mathbf{X})\leq t}|X_i\right]$.

This equation can be rewritten as
\begin{equation}\label{equation:cramer-von-mises}
   CVM_i = \frac{\int \mbox{Var}(\E\left[\1_{f(\mathbf{X})\leq t}|X_i\right])d\mu(t)}{\int \mbox{Var}(\1_{f(\mathbf{X})\leq t}) d\mu(t)}.
\end{equation}

As before, these indices extend to the multivariate case. Simple estimators have been proposed \cite{chatterjee2020new, gamboa2020global}, and are based on permutations and rankings. 

\begin{remark}\label{rmk:cvm_sobol}
    As mentioned earlier, Sobol’ indices quantify correlations and second-order moments but do not take into account information about the distribution of the outcome. However, note the similarity between the definition of the Cramér-von-Mises index and the classical Sobol’ index, especially if we rewrite  Equation \eqref{equation:cramer-von-mises} as:
  \begin{equation}
      CVM_i = \int Sob_i(\1_{f(\mathbf{X})\leq t})\frac{\mbox{Var}(\1_{f(\mathbf{X})\leq t})}{\int \mbox{Var}(\1_{f(\mathbf{X})\leq t}) d\mu(t)}d\mu(t).
  \end{equation} 
  Cramér-von-Mises can be seen as an adaptive Sobol’ index that emphasizes the regions where the cumulative distribution of the outcome is highly changing, as more information can be obtained in these areas. This enable to capture information about the distribution of the outcome instead of moment-related information.
\end{remark}

\subsubsection{Extension of the Cramér-von-Mises indices}

Classical Cramér-von-Mises indices suffer from the same limitation as Sobol’ indices as they are tailored for independent inputs. A natural extension is to create new indices to handle the case of dependent inputs. We propose an extension of the Cramér-von-Mises indices, inspired by the ideas of the extended Sobol’ indices and by the works of \cite{azadkia2019simple}. This new set of indices will capture the influence of a feature independently of the rest of the features.

\begin{definition} For every input $i$, we define the independent Cramér-von-Mises indices as:
\begin{equation}\label{eq:extended_cramer}
\begin{split}
    CVM^{ind}_i  & := \frac{\int \E (\emph{Var}(\1_{f(\mathbf{X})\leq t}|X_{\sim i}))d\mu(t)}{\int \emph{Var}(\1_{f(\mathbf{X})\leq t}) d\mu(t)}\\
\end{split}
\end{equation}
\end{definition}

This extension enables to compare the influence of a feature on the output of an algorithm without its dependencies with other features.

\begin{remark}
  This independent Cramér-von-Mises index can be seen as an extension of the $SobT^{ind}$ index.
\end{remark}
This remark is similar to Remark \ref{rmk:cvm_sobol}. From the independent Total Sobol index shown in \eqref{eq:defSobol’Tind}, by changing the output function as a threshold of the real algorithm and taking the mean along all the possible thresholds, we obtain the independent Cramér-von-Mises index. This index can also be seen as an adaptive form of the $SobT^{ind}$ index.

Estimation of these indices is given in Appendix \ref{Appendix:estimcvm} by the mean of estimates $\widehat{CVM}_i$.  Similarly to Theorem \ref{theo:TCL}, we have the following theorem.
\begin{theorem} If we denote by $N$ the number of observations used to compute $\widehat{CVM}_i$, then the sequence $\sqrt{N}\left(CVM_i - \widehat{CVM}_i\right)$ converges towards the centered Gaussian law with a limiting variance $\xi^2$ whose explicit expression can be found in the proof.
\end{theorem}
The proof of this theorem can be found in \cite{gamboa2018sensitivity}. Note that new estimation procedures can be efficient with little data, as mentioned in \cite{gamboa2020global}, which will be helpful for measuring intersectional fairness in the following Section.

\section{Fairness}\label{section:fairness}
\subsection{Sensitivity Indices as Fairness measures}
In this section, we provide a probabilistic framework to unify  various definitions of  Fairness for Group of individual  as Global Sensitivity Indices. Fairness amounts to quantify the dependencies between a sensitive feature $S$ and functions of the outcome $f(X)$ and of the realisation of the variable of interest $Y$. Several measures of fairness corresponding to different definitions of fairness have been proposed in the machine learning literature.  However, all these definitions boil back to a quantification of the mathematical propositions "$f(X) \ind S$" or "$f(X) \ind S |Y$".

For instance, the two main common definitions of fairness are the following
\begin{itemize}
    \item  \textit{Statistical Parity}, see for instance in~\cite{dwork2012fairness}, requires that the algorithm $f$, predicting a target $Y$, has similar outputs for all the values of $S$ in the sense that the distribution of the output is independent from the sensitive variable $S$, namely  $f(\mathbf{X}) \ind S $. In the binary classification case, it is defined as  $\P(f(\mathbf{X}) = 1 |S ) = \P(f(\mathbf{X}) = 1)$ for general $S$, continuous or discrete.
    \item \textit{Equality of odds} looks for the independence between the error of the algorithm and the protected variable, i.e implying here conditional independence, i.e   $f(\mathbf{X}) \ind S | Y$. This  condition is equivalent in the binary case to $\P(f(\mathbf{X}) = 1 |Y = i, S) = \P(f(\mathbf{X}) = 1 |Y = i),$ for $ i = 0,1.$
\end{itemize}

Previous notions of fairness are quantified using  a \textit{Fairness measure} $\Lambda$ and a function $\Phi(Y,\mathbf{X})$ such that $\Lambda(\Phi(Y,\mathbf{X}),S) = 0$ in the case of perfect fairness while the constraint is relaxed into $\Lambda(\Phi(Y,\mathbf{X}),S) \leq \varepsilon$ , for a small $\varepsilon$, leading to the notion  of approximate fairness. The following definition provides a general framework to define fairness measures. GSA measures as defined in~\ref{section:GSA} or described in~\cite{daveiga:hal-01128666, iooss2015review} are suitable indicators  to quantify  fairness as follows  and these definitions can be extended to continuous predictors and continuous $Y$.
\begin{definition}\label{def:link}
Let $\Phi$ be a function of the features $\mathbf{X}$ and of $Y$. We define a GSA measure for a function $\Phi$ and a random variable $Z$  as a $\Gamma(.,.)$ such that $\Gamma(\Phi(Y,\mathbf{X}),Z)$ is equal to $0$ if $\Phi(Y,\mathbf{X})$ is independent of $Z$ and is equal to $1$ if $\Phi(Y,\mathbf{X})$ is a function of $Z$. Then, $\Gamma$ induces a GSA-Fairness measure defined as  $\Lambda(\Phi(Y,\mathbf{X}),S) = \Gamma(\Phi(Y,\mathbf{X}),S)$.
\end{definition}

The following examples provide a GSA formulation for most of classical fairness definitions using Sobol’ and Cramér-von-Mises indices.

\begin{example}[\textit{Statistical Parity}]\label{ex:SP}
The so-called \textit{Statistical Parity} fairness is achieved by taking $\Lambda(\Phi(Y,\mathbf{X}),S)) = \mbox{Var}(\E[f(\mathbf{X})|S])$. This corresponds to the GSA measure $Sob_S(f(\mathbf{X}))$. If $f$ is a classifier with value in $\{0,1\}$, we recover for a binary $S$ the classical definition of \textit{Disparate Impact},$\P(f(X)=1|S=1) = \P(f(X)=1|S=0)$, see \cite{gordaliza2019obtaining}.
\end{example}

\begin{example}[\textit{Avoiding Disparate Treatment}]
The so-called \textit{Avoiding Disparate Treatment} fairness is achieved by taking $\Lambda(\Phi(Y,\mathbf{X}),S)) = \E[\mbox{Var}(f(\mathbf{X})|X)] $. This corresponds to the GSA measure $SobT_S(f(\mathbf{X}))$. Similarly, for a binary classifier, we recover the classical definition.
\end{example}

\begin{example}[\textit{Equality of Odds}]
The so-called \textit{Equality of Odds} fairness is achieved by taking $\Lambda(\Phi(Y,\mathbf{X}),S)) = \E[\mbox{Var}(\E[f(\mathbf{X})|S,Y]|Y)] $. This corresponds to the GSA measure $CVM^{ind}(f(\mathbf{X}), S|Y)$. Similarly, for a binary classifier, we recover the classical definition.
\end{example}

\begin{example}[\textit{Avoiding Disparate Mistreatment}]\label{ex:ADM}
The so-called \textit{Avoiding Disparate Mistreatment} fairness is achieved by taking $\Lambda(\Phi(Y,\mathbf{X}),S)) = \mbox{Var}(\E[\ell(f(\mathbf{X}),Y)|S]) $ with $\ell$ a loss function. This corresponds to the GSA measure $Sob_S(\ell(f(\mathbf{X}),Y))$. Similarly, for a binary classifier, we recover the classical definition.
\end{example}

 Among well known fairness measures, we point out that we immediately recover two main fairness measures used in the fair learning literature -- namely \textit{Statistical Parity} and \textit{Equality of Odds}. GSA measures can be computed for different function $\Phi$ and highlight either the behaviour of the algorithm, $\Phi(Y,\mathbf{X}) = f(\mathbf{X})$, or its performance, $\Phi(Y,\mathbf{X}) = \ell(Y,f(\mathbf{X}))$ for a given loss $\ell$. This can lead to different GSA-Fairness definitions from a same GSA measure, see Examples \ref{ex:SP} and \ref{ex:ADM}.

\begin{example}
  Recent work in Fairness literature exposed various definitions and measures to quantify influence of a sensitive feature, beyond classical notions. For instance, \cite{fairness_shap} uses Shapley values, \cite{li2019kernel} uses HSIC measures, \cite{ghassami2018fairness} uses Mutual Information, so on and so forth. All these measures have been extensively studied in GSA literature, as mentioned in previous Section, and these frameworks are included in ours. 
\end{example}

In Table \ref{tab:Fairness_GSA}, we summarize the different indices associated to classical studied fairness definitions shown in previous Examples. By considering these fairness definitions as GSA measures, we can explain fairness in terms of simple effects presented in previous section, along with limitations of those definitions. For instance, \textit{Statistical Parity} corresponds to the classical Sobol’ index. The nullity of this index implies no direct influence of sensitive variables on the outcome, but can be limited as sensitive variables may have joint effects with other variables not captured by this metric. Therefore, \textit{Statistical Parity} will lack in this regard. On the contrary, since \textit{Avoiding Disparate Treatment} corresponds to Total Sobol’ indices, this definition of fairness captures every possible influence of the sensitive feature on the outcome.
\begin{table*}[t]
\caption{Common fairness definitions and associated GSA measures}
\vskip 0.15in
\centering
\begin{sc}
\begin{tabular}{ |c|c| }
 \hline
 Fairness definition & GSA measure associated  \\
 \hline
 Statistical Parity & $\mbox{Var}(\E[f(\mathbf{X})|S]) \to Sob_S(f(\mathbf{X}))$  \\
 \hline
 Avoiding Disparate Treatment & $\E[\mbox{Var}(f(\mathbf{X})|X)] \to SobT_S(f(\mathbf{X}))$  \\
 \hline
 Equality of odds & $\E[\mbox{Var}(\E[f(\mathbf{X})|S,Y]|Y)] \to CVM^{ind}(f(\mathbf{X}), S|Y)$  \\
 \hline
 Avoiding Disparate Mistreatment & $\mbox{Var}(\E[\ell(f(\mathbf{X}),Y)|S]) \to Sob_S(\ell(f(\mathbf{X}),Y))$  \\
 \hline
\end{tabular}
\end{sc}
\label{tab:Fairness_GSA}
\end{table*}

\begin{remark}
Note that many fairness measures are defined using discrete or binary sensitive variable. The GSA framework enables to handle continuous variables without additional difficulties. Moreover using kernel methods, GSA indices can be defined for a larger and more "exotic" variety  of variables such as graphs or trees, for instance. In particular  HSIC (see in \cite{daveiga:hal-01128666, berlinet2004collection, gretton2005kernel, smola2007hilbert, meynaoui2019new}) is a kernel-based GSA measure that has been used in fairness.
\end{remark}

\subsection{Consequences of seeing Fairness with  Global Sensitivity Analysis optics}

In this subsection, we enumerate various consequences of studying Fairness with this probabilistic framework coming from the GSA literature. 

\begin{enumerate}[(i)]
  \item \textbf{Modularity of fairness indicators}
  Numerous metrics have been proposed in GSA literature to quantify the influence of a feature on the outcome of an algorithm. We already mentioned several of them so far. This diversity enables choices in the quantified fairness since every choice of GSA measure induces a Fairness definition. We presented in previous subsection a concrete example with Sobol’ indices, namely between \textit{Disparate Impact} and \textit{Avoiding Disparate Treatment}. Another example would be the use of kernels in HSIC-based indices, as exposed for instance in \cite{li2019kernel}. By selecting various kernels, specific characteristics associated with fairness can be targeted.
  \item \textbf{Perfect and Approximate fairness}
  GSA has been especially created to quantify \textit{quasi} independence between variables. Merging GSA and Fairness gives a formal framework to the notion of approximate fairness and computationally justify the use of GSA codes to measure and quantify fairness. Additionally, as mentioned in previous section, GSA literature includes statistical tests for independence between input variables and outcomes, along with confidence intervals. Therefore, it is possible to compute them in order to test whether perfect fairness or approximate fairness is obtained. Moreover, this enables the possibility of auditing algorithms.
  \item \textbf{Choice of the target}
  The framework presented earlier works for quantifying the influence of a sensitive feature on the outcome of a predictor but also any function of the predictor and of the input variables. This includes the loss of a predictor against a target. The ambivalence of this framework allows links to be made between various fairness definitions. For example, \textit{Disparate Impact} and \textit{Avoiding Disparate Mistreatment} are the same fairness but applied either to the predictor or to the loss of the predictor against a real target. In the first case, we want the algorithm to be independent of the sensitive feature; while in the second case, we want the errors of the predictor to be independent of the sensitive feature. Moreover, it allows for extension of fairness definitions to cases where an algorithm can be biased, as long as it does not make a mistake. 
  \item \textbf{Second-level Global Sensitivity Analysis}
  Recent works in GSA take into account the uncertainty of the distribution of the inputs of an algorithm, see \cite{meynaoui2019new}. These tools can help in a fairness framework, especially when the distribution of sensitive features is unknown and unreachable. This will be more deeply studied in future papers.
\end{enumerate}

\subsection{Applications to Causal Models}

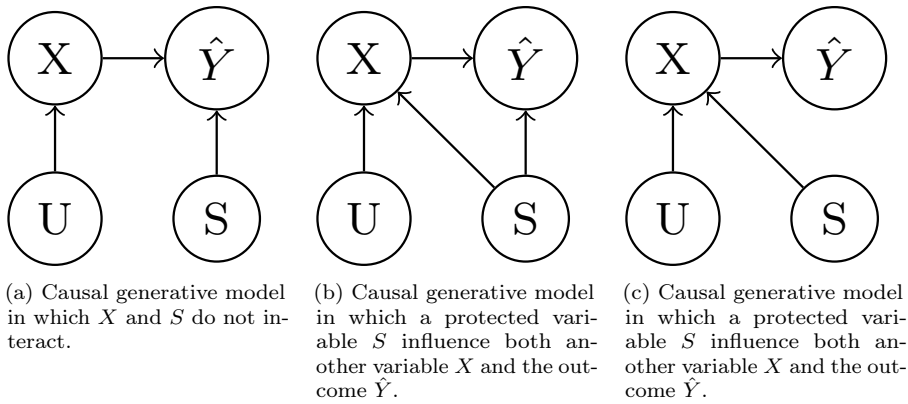
\begin{figure}
\hfill
    \begin{subfigure}[t]{0.3\linewidth}
        \centering
        \resizebox{\linewidth}{!}{\begin{tikzpicture}
          \node[shape=circle,draw=black] (S) at (1,0) {S};
          \node[shape=circle,draw=black] (X) at (0,1) {X};
          \node[shape=circle,draw=black] (U) at (0,0) {U};
          \node[shape=circle,draw=black] (Y) at (1,1) {$\hat{Y}$};

          \path [->] (U) edge node[left] {} (X);
          \path [->] (S) edge node[left] {} (Y);
          \path [->] (X) edge node[left] {} (Y);
        \end{tikzpicture}
}
        \caption{Causal generative model in which $X$ and $S$ do not interact.}
        \label{figure:expe_buhlmann_1}
    \end{subfigure}
    \hfill
    \begin{subfigure}[t]{0.3\linewidth}
        \centering
        \resizebox{\linewidth}{!}{\begin{tikzpicture}
          \node[shape=circle,draw=black] (S) at (1,0) {S};
          \node[shape=circle,draw=black] (X) at (0,1) {X};
          \node[shape=circle,draw=black] (U) at (0,0) {U};
          \node[shape=circle,draw=black] (Y) at (1,1) {$\hat{Y}$};

          \path [->] (S) edge node[left] {} (X);
          \path [->] (U) edge node[left] {} (X);
          \path [->] (S) edge node[left] {} (Y);
          \path [->] (X) edge node[left] {} (Y);
        \end{tikzpicture}
}
        \caption{Causal generative model in which a protected variable $S$ influence both another variable $X$ and the outcome $\hat{Y}$.}
        \label{figure:expe_buhlmann_2}
    \end{subfigure}
    \hfill
    \begin{subfigure}[t]{0.3\linewidth}
        \centering
        \resizebox{\linewidth}{!}{\begin{tikzpicture}[ auto]
          \node[shape=circle,draw=black] (S) at (1,0) {S};
          \node[shape=circle,draw=black] (X) at (0,1) {X};
          \node[shape=circle,draw=black] (U) at (0,0) {U};
          \node[shape=circle,draw=black] (Y) at (1,1) {$\hat{Y}$};

          \path [->] (S) edge node[left] {} (X);
          \path [->] (U) edge node[left] {} (X);
          \path [->] (X) edge node[left] {} (Y);
        \end{tikzpicture}
}
        \caption{Causal generative model in which a protected variable $S$ influence both another variable $X$ and the outcome $\hat{Y}$.}
        \label{figure:expe_buhlmann_3}
    \end{subfigure}
    \hfill
    \caption{Examples of representation of causal models with directed acyclic graphs.}
    \label{figure:expe_buhlmann}
\end{figure}

Quantifying fairness using measures is a first step to understand bias in Machine Learning. Yet, causality enables to understand the true reasons of discrimination, as it is often related to the causal effect of a variable. The relations between variables describing causality are often modeled using a Directed Acyclic Graph (DAG). We refer to \cite{pearl2009causality, bongers2020foundations}.  

In this subsection, we show how to address causal notions of fairness using the GSA framework, illustrated by a synthetic and a social example. We show that information gained thanks to Sobol’ indices allow to learn some characteristic about the causal model.

We tackle the problem of predicting $Y$ by $\hat{Y}$ knowing $(X, S)$ while the non-sensitive variables are influenced by a non-observed exogeneous variable $U$. This is modeled by the following equations:

\[
  X = \phi(U,S)\quad
  \hat{Y} = \psi(X,S),\]
where $\phi$ and $\psi$ are some unknown functions. These equations are a consequence of the unique solvability of acyclic models \cite{bongers2020foundations} and are illustrated in the various DAGs of Figure \ref{figure:expe_buhlmann}.

In many practical cases, the causal graph is unknown and we need indices to quantify causality.  In the following, we are not interested in the complete knowledge of the graph -- which is a NP-hard problem -- but only in the existence of paths from $S$ to $Y$.\vskip .1in Actually, GSA can quantify causal influence following DAG structure, and different GSA indices will correspond to different paths from $S$ to $Y$.   Different type of relationships can be measured in particular with the Total Sobol and the Total Independent Sobol indices to quantify  either the presence of a path from $S$ directly to $Y$ or a path from $S$ to another variable $X$ that influences itself the predictor $Y$. We call this latter effect a "bouncing effect" since $Y$ is influential only through a mediator.  \\
The following proposition explains how specific Sobol indices can be used to detect the presence of causal links between the sensitive variable and the outcome of the algorithm.

 \begin{proposition}[Quantifying Causality with Sobol Index] \label{prop:causal}$ $ \\
\begin{itemize}
\item The condition $SobT_S = 0$ implies that every path from $S$ to $Y$ is non-existent, that is $S$ and $Y$ belong to two different connected component of the causal graph. 
\item The condition $SobT_S^{ind} = 0$ implies that the direct path from $S$ to $Y$ is non-existent, that is the absence of direct edge between $S$ and $Y$ in the causal graph.
\end{itemize}
\end{proposition}

Hence, using GSA, we can infer the absence of causal link between sensitive features and outcomes of algorithm without knowing the structure of the DAG. Note that, while Sobol’ indices are correlation-based, this is not an issue in quantifying causality for fairness, as the sensitive features are usually supposed to be roots of the DAG \cite{bongers2020foundations,de2021counterfactual}.

\begin{example}[Causal graphs \cite{rothenhausler2018anchor}]\label{ex:buhlmann}
In this example, we specify three causal models and illustrate the previous proposition.
 
 In Graph \ref{figure:expe_buhlmann_1}, $S$ is directly influent on the outcome $\hat{Y}$. There is no interaction between $S$ and $X$. This happens when $S$ and $X$ are independent for instance. In such a case, Sobol’ indices and independent Sobol’ indices are the same, as mentioned in Remark \ref{remark:sobol=sobolind}. The equality $SobT_S = SobT_S^{ind}$ ensures the absence of "bouncing effect" for the sensitive variable $S$.
 
 In Graph \ref{figure:expe_buhlmann_2}, we have no information about the influence of $S$ on the outcome.
 
 In Graph \ref{figure:expe_buhlmann_3}, $S$ has no direct influence on the outcome, therefore $SobT_S^{ind} = 0$. This variable can still be influent on the outcome since it may modify other variables of interest. In this case, $X$ is a mediator variable through which the sensitive feature will influence the outcome with a "bouncing effect".  A model describing this kind of DAG in a fairness framework is the "College admissions" case, explained below.
\end{example}

\begin{example}[College admissions]\label{ex:admissions}
  This example focus on college admissions process. Consider $S$ to be the gender, $X$ the choice of department, $U$ the test score and $\hat{Y}$ the admission decision. The gender should not directly influence any admission decision $\hat{Y}$, but different genders may apply to departments represented by the variable $X$ at different rates, and some departments may be more competitive than others. Gender may influence the admission outcome through the choice of department but not directly.  In a fair world, the causal model for the admission can be modeled by a DAG  without direct edge from $S$ to $\hat{Y}$. Conversely, in an unfair world, decisions can be influenced directly by the sensitive feature $S$ -- hence the existence of a direct edge between $S$ and $\hat{Y}$. This issue on
  \textit{unresolved discrimination} is tackled in~\cite{kilbertus2017avoiding, frye2020asymmetric}.
\end{example}

\begin{table}[t]
\caption{Sobol’ indices: what is taken into account and what is not.}
\vskip 0.15in
\centering
\begin{sc}
\begin{tabular}{ |c|c|c|  }
 \hline
 \multicolumn{3}{|c|}{Sobol’ indices} \\
 \hline
 &  Correlation between & Joint\\
 &  variables & contributions\\
 \hline
 $Sob_i$ & \cmark & \xmark\\
 \hline
 $SobT_i$ & \cmark & \cmark\\
 \hline
 $Sob^{ind}_i$ & \xmark &\xmark\\
 \hline
 $SobT^{ind}_i$ &\xmark & \cmark\\
 \hline
\end{tabular}
\end{sc}

\end{table}

\subsection{Quantifying intersectional (un)fairness with GSA index}

Most of fairness results are stated in the case where there is only one sensitive variable. Yet in many cases, the bias and the resulting possible discrimination are the result of multiple  sensitive variables. This situation is known as intersectionality,  when the level of discrimination of an intersection of several minority groups is worse than the discrimination present in each group as presented in \cite{crenshaw1989demarginalizing}. Some recent works provide extensions of fairness measures to take into account the bias amplification due to intersectionality. We refer for instance to \cite{morina2019auditing}  or \cite{foulds2020intersectional}. However, quantifying this worst case scenario cannot be achieved using standard fairness measures. The GSA framework allows for controlling the influence of a set of variables and as such can naturally address intersectional notions of fairness. \\


Intersectional fairness is obtained when multiple sensitive variables (for instance $S_1$ and $S_2$ in the most simple case) do not have any joint influence on the output of the algorithm. We propose a definition of intersectional fairness using GSA indices.

\begin{definition}
Let $S_1, S_2, \cdots, S_m$ be sensitive features.
It is said that an algorithm output is intersectionaly fair if $\Gamma(\Phi(X, S_1, \cdots, S_m); (S_1, \cdots, S_m)) = 0$. This constraint can be relaxed to $\Gamma(\Phi(X, S_1, \cdots, S_m); (S_1, \cdots, S_m)) \leq \varepsilon$ with $\varepsilon$ small for approximate intersectionality fairness.
\end{definition}

Consider two independent protected features $S_1$ and $S_2$ (i.e gender and ethnicity). Depending on the chosen definition of fairness, there are situation where fairness is obtained with respect to $S_1$, with respect to $S_2$ but where the combined effect of $(S_1, S_2)$ is not taken into account.  For instance, let $Y=S_1\times S_2$. In this toy-case, the Disparate Impact of $S_1$, as well as the Disparate Impact of $S_2$, is equal to $1$ while the Disparate Impact of $(S_1, S_2)$ is equal to $0$. This can be readily seen thanks to the link between fairness and GSA as the Sobol’ indices for $S_1$ and for $S_2$ are null while the Sobol’ index for the couple $(S_1, S_2)$ is maximal.

\begin{proposition}\label{prop:intersect1}
Let $(S_1, S_2, \cdots, S_m)$ be sensitive features. To be fair in the sense of \textit{Disparate Impact} for $S_1$ and to be fair in the sense of \textit{Disparate Impact} for $S_2$ does not quantify any intersectional fairness in the sense of the \textit{Disparate Impact}.
\end{proposition}

However, if we take again the same toy-case but look at the Total Sobol’ indices, we see that $SobT_{S_1} = 0$ implies that $SobT_{(S_1, S_2)} =0$.

\begin{proposition}\label{prop:intersect2}
Let $(S_1, S_2, \cdots, S_m)$ be sensitive features. To be fair in the sense of \textit{Avoiding Disparate Treatment} for $S_1$ implies intersectional fairness for any intersection where $S_1$ appears.
\end{proposition}

\begin{remark}
Intersectional fairness is different than classical fairness. Classical fairness only pays attention to the influence of a single sensitive feature on the outcome while intersectional fairness is quantifying only the influence due to interactions between sensitive features. In applications, the goal is usually to have both classical and intersectional fairness. A single fairness definition that covers these two characteristics can be hard to find or too restrictive to readily use. For instance, among Sobol’ indices, only the Total Sobol’  index induces both a classical and intersectional fairness.
\end{remark}

\section{Experiments}\label{section:experiments}

\begin{table}[t]
    \caption{Synthetic experiments based on causal DAGs -- Figure \ref{figure:expe_buhlmann}}. 
    \vskip 0.15in
    \centering
    \begin{sc}
    \begin{tabular}{ |c|c|c|c|c| }
         \hline
         & $Sob$ & $SobT$ & $Sob^{ind}$ & $SobT^{ind}$ \\
         \hline
         \hline
         \multicolumn{5}{|c|}{$Y = 2\times X $}\\
         \hline
         X & \textbf{1.00} (0.99 - \textbf{1.00} - 1.00) & \textbf{1.00} (0.99 - \textbf{1.00} - 1.00) & \textbf{0.75} (0.74 - \textbf{0.75} - 0.76) & \textbf{0.75} (0.74 - \textbf{0.75} - 0.76) \\
         \hline
         S & \textbf{0.24} (0.24 - \textbf{0.25} - 0.26) & \textbf{0.25} (0.24 - \textbf{0.25} - 0.26) & \textbf{0.00} (0.00 - \textbf{0.00} - 0.01) & \textbf{0.00} (0.00 - \textbf{0.00} - 0.01)\\
        \hline
        \hline
        \multicolumn{5}{|c|}{$Y = 0.7\times X + 0.3\times S$}\\
         \hline
        X & \textbf{0.91} (0.89 - \textbf{0.91} - 0.93) & \textbf{0.92} (0.89 - \textbf{0.91} - 0.94) & \textbf{0.51} (0.46 - \textbf{0.48} - 0.52) & \textbf{0.52} (0.46 - \textbf{0.47} - 0.54)\\
         \hline
          S & \textbf{0.52} (0.48 - \textbf{0.53} - 0.55) & \textbf{0.54} (0.48 - \textbf{0.53} - 0.55) & \textbf{0.07} (0.05 - \textbf{0.09} - 0.11) & \textbf{0.09} (0.06 - \textbf{0.09} - 0.12)\\
         \hline
         \hline
         \multicolumn{5}{|c|}{$Y = 0.7\times X + 0.3\times S$}\\
         \hline
         X & \textbf{0.78} (0.78 - \textbf{0.84} - 0.85) & \textbf{0.84} (0.80 - \textbf{0.84} - 0.86) & \textbf{0.81} (0.78 - \textbf{0.84} - 0.85) & \textbf{0.82} (0.80 - \textbf{0.84} - 0.86)\\
         \hline
          S & \textbf{0.13} (0.12 - \textbf{0.16} - 0.17) & \textbf{0.17} (0.15 - \textbf{0.16} - 0.18) & \textbf{0.14} (0.12 - \textbf{0.16} - 0.17) & \textbf{0.15} (0.13 - \textbf{0.16} - 0.18)\\
         \hline
    \end{tabular}
    \caption*{Legend: Values format is "\textbf{experimental value} (lower bound of 95\% confidence interval - \textbf{theoretical value} - upper bound of 95\% confidence interval)". }
\end{sc}
\label{tab:Experiments}
\end{table}

\subsection{Synthetic experiments}
In this subsection, we focus on the computation of complete Sobol’ indices in a synthetic framework. We design three experiments, modeled after the causal generative models shown in Figure \ref{figure:expe_buhlmann}. For simplicity, we consider a Gaussian model. In each experiment $j,j \in \{1, 2, 3\}$, $(X,S,U)$ are random variables drawn from a Gaussian distribution with covariance matrix $C_j$, where
\begin{equation*}
    C_1 = C_2 = \begin{pmatrix}
    1 & 0.5 & 0.5\\
    0.5 & 1 & 0 \\
    0.5 & 0 & 1 \\
    \end{pmatrix}, C_3 = \begin{pmatrix}
    1 & 0 & 0.5\\
    0 & 1 & 0 \\
    0.5 & 0 & 1 \\
    \end{pmatrix}.
\end{equation*}
The random variable $U$ is unobserved in this case and therefore does not have Sobol’ indices. Its purpose is to simulate exogenous variables that modify the features in $X$. The target $Y_j$, described in the Table \ref{tab:Experiments} for each of the experiments, is equal to
\begin{align*}
    Y_1 = & 2 \times X,\\
    Y_2 = Y_3 = & 0.7 \times X + 0.3 \times S.\\
\end{align*}
The first experiment shows the difference between independent and non-independent Sobol’ indices. The outcome is entirely determined by a single variable $X$ and therefore, $Sob_X = 1$. However, $X$ is intrinsically linked with a sensitive feature because of the covariance matrix, so that $Sob_X^{ind} \not = 0$. This is a concrete example where \textit{Statistical parity} is not obtained for $S$ but \textit{unresolved discrimination} mentioned in Example \ref{ex:admissions} is obtained, since $S$ is influential only through $X$.

The second experiment adds a direct path from the variable $S$ to the outcome $Y$. Since $Y$ can be factorized as an effect from $X$ and an effect of $S$, we still have $Sob_X = SobT_X$ and $Sob_X^{ind} = SobT^{ind}_X$. However, in this case, $X$ is no longer enough to fully explain the outcome, so that $Sob_X \not = 1$. $Sob_S^{ind}$ quantify the influence of this direct path from $S$ to $Y$. Note that the difference between $Sob_S$ and $Sob^{ind}_S$ quantify the influence of the path from $S$ to $Y$ through the intermediary variable $X$.

In the third experiment, $S$ and $X$ are independent and $S$ can only influence the outcome directly. This is the framework of classical Global Sensitivity Analysis. In this case, non-independent and independent Sobol’ indices are equal, as mentioned in Remark \ref{remark:sobol=sobolind}

\begin{figure}
    \centering
    \includegraphics[scale=.5]{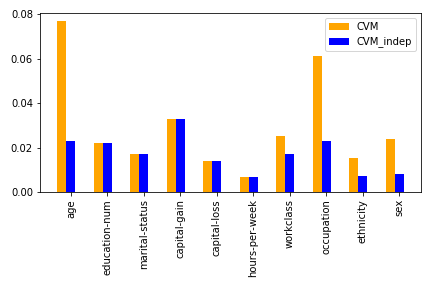} 
    \caption{Cramér-von-Mises and independent Cramér-von-Mises indices for the Adult dataset.}
    \label{exp:adult}
\end{figure}
\subsection{Real data sets}
In this section, we focus on the implementation of Cramér-von-Mises indices on two real-life datasets: the Adult dataset \cite{Dua:2019} and the COMPAS dataset.

\subsubsection{Adult dataset}
The adult  dataset   consists in   14   attributes  for 48,842 individuals. The class label corresponds to the annual income (below/above 50.000 $k\$$). We  study the effect of  different attributes. The results for a classifier obtained for an algorithm built using an Extreme Gradient Boosting Procedure are shown in Figure~\ref{exp:adult}. We used the same pre-process as \cite{besse2020survey} for the choice of variables.\\
If we look at the independent Cramér-von-Mises, we quantify the direct influence of a variable . We recover the influent indicators -- "capital gain", "education-number", "age", "occupation"... -- given by other studies \cite{frye2020asymmetric, besse2020survey}. \\
The joint influences on the outcome of other variables is also measured using GSA indices. Variables for which independent and classical Cramér-von-Mises indices are the same have no "bouncing" influence. Otherwise, the gap between these two indices quantify this specific effect. For example, the variable "age"  correlates with most of the other variables such as "education-number" or "marital-status" for instance. Because of this, most of its influence is through "bouncing effects" and the gap between its two indices (i.e "$CVM$" and "$CVM_indep$") is larger than for any other feature. 
The variable "sex" also plays an important role through its "bouncing" effect. We can see this through the difference between the classical and the independent index associated with this feature. This explains why removing the variable "sex" is not enough to obtain a fair predictor since it influences other variables that affect the prediction. We recover the results obtained by several studies that point out the bias created by the "sex" variable.\\ Note that race may have led to unbalanced decisions as well. Yet, the Cramér-von-Mises index is lower than the one for the "sex" variable, which explains why the discrimination is lower than the one created by the sex, as emphasized by the study of the Disparate Impact  which is in a 95\% confidence interval of $[0.34, 0.37]$ for sex and $[0.54, 0.63]$ for ethnic origin in \cite{besse2020survey}.

\subsubsection{COMPAS dataset}
\begin{figure}
    \centering
    \begin{subfigure}[t]{0.49\linewidth}
        \centering
        \includegraphics[width = \textwidth]{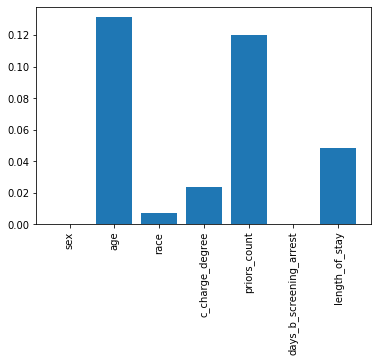}
        \label{figure:compas}
        \caption{Cramér-von-Mises indices computed for the COMPAS decile score.}
    \end{subfigure}
    \hfill
    \begin{subfigure}[t]{0.49\linewidth}
        \centering
        \includegraphics[width = \textwidth]{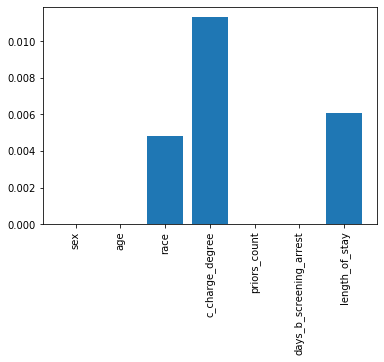}
        \label{figure:compasloss}
        \caption{Cramér-von-Mises indices computed on the loss between COMPAS output and real case of recidivism after two years.}
    \end{subfigure}
    \caption{Cramér-von-Mises indices for the COMPAS dataset.}
    \label{exp:compas}
\end{figure}

The so-called COMPAS  dataset, gathered by ProPublica described for instance in in \cite{washington2018argue} , contains information about the recidivism risk predicted by the COMPAS tool, as well as the ground truth recidivism rates, for 7214 defendants. The  COMPAS risk score,  between $1$ and $10$ ($1$ being a low chance of recidivism and $10$ a high chance of recidivism),  is obtained by an algorithm using all other variables used to compute it, and is used to forecast  whether the defendant will  reoffend or not.
We analysed this dataset with Cramér-von-Mises indices in order to quantify fairness exhibited by the COMPAS algorithm. The results are shown in Figure \ref{exp:compas}.

First, every independent index is null, which means that the COMPAS algorithm does not rely on a single variable to predict recidivism. Also, gender and ethnicity are virtually not used by the algorithm, opposed to the variables "age" or "$priors\_count$" (the number of previous crimes). Hence as expected, the algorithm appears to be fair. However, when comparing the accuracy  of the predictions of the algorithm with  real-life two-year recidivism, the "race" variable is found to be influential. Hence we show that the indices we propose recover the bias denounced by Propublica with an algorithm that, despite fair predictions,  shows a behavior that favors a part of the population based on the race variable.

\section{Conclusion}\label{section:conclusion}
We recalled classical notions both for the Global Sensitivity Analysis and the Fairness literature. We presented new Global Sensitivity Analysis tools by the mean of extended Cramér-von-Mises indices, as well as proved asymptotic normality for the extended Sobol’ indices. These sets of indices allow for uncertainty analysis for non-independent inputs, which is a classical situation in real-life data but not often studied in the literature. Concurrently, we link Global Sensitivity Analysis to Fairness in an unified probabilistic framework in which a choice of fairness is equivalent to a choice of GSA measure. We showed that GSA measures are natural tools for both the definition and comprehension of Fairness. Such a link between these two fields offers practitioners customized techniques for solving a wide array of fairness modeling problems.

\begin{acknowledgements}
Research partially supported by the AI Interdisciplinary Institute ANITI, which is funded by the French “Investing for the Future – PIA3” program under the Grant agreement ANR-19-PI3A-0004.
\end{acknowledgements}

%
\section*{Conflict of interest}
The authors declare that they have no conflict of interest.

\bibliographystyle{spmpsci}      
\bibliography{biblio.bib}   

%
%

\appendix

\section{Lévy-Rosemblatt theorem and associated mappings}\label{Appendix:rosemblatt}
The aim of the Lévy-Rosenblatt transform is to find a transport map between the correlated $\mathbf{X}$ and independent uniform variables $\mathbf{U} \in \R^p$. From now, we assume the distribution of $\mathbf{X}$ to be absolutely continuous.
\begin{theorem}[Lévy-Rosemblatt theorem,\cite{levy1954theorie, rosenblatt1952}]: \label{levy} there is a bijection (denoted "RT" for Rosemblatt transform) between $p(\mathbf{X})$ and $p$ independent uniform random variables
\begin{equation}\label{Rosenblatt_transform}
        (X_i, (X_{i+1}|X_i), \cdots, (X_{i-1}|X_{\sim (i-1)})) \sim p_\mathbf{X} \xrightarrow{RT} (U^i_1,\cdots,U^i_p) \sim \mathcal{U}^p(0,1).
\end{equation}
\end{theorem}

\begin{example}
  In the following, we will always be interested in two groups of variables: the sensitive variable $X_i$ and the rest of the variables $X_{\sim i}$. Therefore, it may help to understand the special case where $\mathbf{X} = (X_1,X_2)$ since it encapsules all the difficulty. In this case, we have two different ways to decompose $p_{\mathbf{X}}$.
  \begin{enumerate}[(i)]
    \item If we decompose $p_{\mathbf{X}}$ as $p_{X_1}\times p_{X_2|X_1}$, then we can map this to $(U^1_1, U^1_2)$. With this mapping, we can draw random variables with distributions $p_{X_1}$ and $p_{X_2|X_1}$. For this, we need only to have access to independent uniform random variables and use the inverse Rosenblatt transform. We denote as $F_T$ the cumulative distribution function of the random variable $T$. The inverse Rosenblatt transform is then given by
    \begin{align}
      z_1 &= F^{-1}_{X_1}(u^1_1)\\
      z_2 &= F^{-1}_{X_2|X_1 = x_1}(u^1_2).
    \end{align}
    We first draw a random variable $Z_1$ with distribution $p_{X_1}$ from an uniform random variable by quantile inversion. Now that we have this realisation $z_1$, we have the second distribution $p_{X_2|X_1 = z_1}$. We then draw a random variable $Z_2$ that follows the distribution $p_{X_2|X_1 = z_1}$ and such that the couple $(Z_1,Z_2)$ has the same distribution as $(X_1,X_2)$. This random variable is similar to $X_2$ but does not contain its correlation with $X_1$.
    \item Similarly, if we decompose $p_{\mathbf{X}}$ as $p_{X_2}\times p_{X_1|X_2}$, then we can map this to $(U^2_1, U^2_2)$.
  \end{enumerate}
  Note that the only case where these two mappings are similar is when $X_1$ and $X_2$ are independent. In that case, $p_{X_1} = p_{X_1|X_2}$ and $p_{X_2} = p_{X_2|X_1}$.
\end{example}

Several things need to be said about this transform.
\begin{remark}
    It enables to transform a set of possibly dependent random variables into a set of random variables without any dependencies. Moreover, for one such set of independent variables $\mathbf{U}^i$, there exists a function $g_i$ square integrable such that $f(\mathbf{X}) = g_i(\mathbf{U}^i)$. One way to compute Sobol’ indices for the output $f(\mathbf{X})$ is therefore to use the Hoeffding decomposition of $g_i(\mathbf{U}^i)$.
\end{remark}
\begin{remark}
     In terms of information, $U^i_1$ carries as much information as $X_i$ since $U^i_1 = F_{X_i}(X_i)$. Note that this include the eventual dependency with other variables. This means that the Sobol’ indices of $U^i_1$ will correspond to the Sobol’ indices of $X_i$ as defined in the previous section. Meanwhile, the law of $U^i_n$ is associated with the law of $X_{i-1}|X_{\sim (i-1)}$. This conditional distribution aim to capture all the remaining randomness in $X_{i-1}$ when the intrinsic effects of the others inputs on it has been removed. Therefore, it has all the remaining information in the law of $X_{i-1}$ when the contribution of the other variables are discarded.
\end{remark}

\begin{remark}
    The previous point is the reason why we do not need to consider all $n!$ possible Rosenblatt Transforms of $\mathbf{X}$. Since we are only interested in the information carried by a variable -- with ($X_i$)  -- and by the law of this same variable without its dependencies in the other variables -- with ($X_i|X_{\sim i}$), we are only interested in $U^i_1$ and $U^i_n$, for all $i$. Therefore, we can without loss of generality, consider a cyclic permutation. That being said, if, for numerical reasons, other Rosenblatt transforms are easier to work with, there is no theoretical reasons not to use them.
\end{remark}

In the classic Sobol’ analysis, for an input $Y$, we have two indices that quantify the influence of the considered feature on the output of the algorithm, namely the first order and total indices. Now, thanks to the Lévy-Rosemblatt, we have two different mappings of interest: the mapping from $U^i_1$ to $X_i$ that includes the intrinsic influence of other inputs over this particular input and the mapping from $U^{i+1}_p$ to $X_i|X_{\sim i}$ that excludes these influences and shows the variation induced by this input on its own. These two different mappings will each lead to two indices (the Sobol’ and Total Sobol’ indices of $U^i_1$, and the ones of $U^{i+1}_p$) so every input $X_i$ will be represented by four indices.

\section{Estimates of extended Sobol’ indices}\label{Appendix:estim_mara}
We recall that in the independent Sobol’ framework, for every input $X_k$, we have two different mappings: the mapping from $U^k_1$ to $X_k$ that includes the intrinsic influence of other inputs over this particular input and the mapping from $U^{k+1}_p$ to $X_k|X_{\sim k}$ that excludes these influences and shows the variation of this input on its own. These two different mappings will each lead to two indices (the Sobol indices of $U^k_1$ and the ones of $U^{k+1}_p$) so every input $X_k$ will be represented by four indices, explained in the following subsection.

As seen previously, the four Sobol’ indices for each variable $X_i , i \in \llbracket1,n \rrbracket$  are defined as followed:
\begin{equation}
    Sob_i = \frac{V[\E[g_i(\mathbf{U}^i)|U^i_1]]}{V[g_i(\mathbf{U}^i)]}  = \frac{V[\E[f(\mathbf{X})|X_i]]}{V[f(\mathbf{X})]}
\end{equation}
\begin{equation}
    SobT_i  = \frac{\E[V[g_{i}(\mathbf{U}^{i})|U^{i}_{\sim 1}]]}{V[g_{i}(\mathbf{U}^{i})]} = \frac{\E[V[f(\mathbf{X})|Z_i]]}{V[f(\mathbf{X})]}
\end{equation}
\begin{equation}
    Sob_i^{ind} = \frac{V[\E[g_{i+1}(\mathbf{U}^{i+1})|U^{i+1}_p]]}{V[g_{i+1}(\mathbf{U}^{i+1})]}  = \frac{V[\E[f(\mathbf{X})|Z_i]]}{V[f(\mathbf{X})]}
\end{equation}
\begin{equation}
    SobT_i^{ind}  = \frac{\E[V[g_{i+1}(\mathbf{U}^{i+1})|U^{i+1}_{\sim p}]]}{V[g_{i+1}(\mathbf{U}^{i+1})]} = \frac{\E[V[f(\mathbf{X})|X_{\sim i}]]}{V[f(\mathbf{X})]}
\end{equation}

We recall that these indices use the Rosemblatt transform, a bijection between independent uniforms and the distribution of the features. This bijection can be inverted to generate samples from uniforms. We denote the inverse of the Rosemblatt transform as IRT -- Inverse Rosemblatt Transform. Thanks to the IRT, we can generate four samples:

\begin{align}
\label{equation:monte_carlo_samples_sobol}
\begin{split}
 (u^i_1,\cdots,u^i_p) &\xrightarrow{\emph{IRT}} \mathbf{x} = (x_i,\cdots,x_{i-1}) \sim p(\mathbf{X}) ,
\\
(u^{i\prime}_1,\cdots,u^{i\prime}_p) &\xrightarrow{\emph{IRT}} \mathbf{x}^{\prime}  = (x^\prime_i,\cdots,x^\prime_{i-1}) \sim p(\mathbf{X}) ,
 \\
(u^i_1,u^{i\prime}_2,\cdots,u^{i\prime}_p) &\xrightarrow{\emph{IRT}} \mathbf{x}^i =  (x_i,x^\prime_{i+1}\cdots,x^\prime_{i-1}) \sim p(X_i)p(X_{\sim i}|X_i) ,
 \\
(u^{i\prime}_1,\cdots,u^{i\prime}_{p-1},u^{i}_p) &\xrightarrow{\emph{IRT}} \mathbf{x}^{i-1} = (x^\prime_i,x^\prime_{i+1}\cdots,x_{i-1}) \sim p(X_{\sim i-1})p(X_{i-1}|X_{\sim i-1}) .
\end{split}
\end{align}

Once we obtain, for each $i \in \{1,\cdots,p\}$, the four samples defined above, we can compute the estimators of the Sobol’ and independent Sobol’ indices as follows:

\begin{align}
\begin{split}
    \widehat{Sob_i} &= \frac{\frac{1}{N} \sum^N_{k = 1}f(\mathbf{x}_k) \times \left(f(\mathbf{x}^i_k) - f(\mathbf{x}_k')\right)}{\hat{V}}\\
    \widehat{SobT^{ind}_i} &= \frac{\frac{1}{N} \sum^N_{k = 1} \left(f(\mathbf{x}^{i-1}_k) - f(\mathbf{x}_k')\right)^2}{2\hat{V}}\\
    \widehat{Sob^{ind}_{i-1}} &= \frac{\frac{1}{N} \sum^N_{k = 1}f(\mathbf{x}_k) \times \left(f(\mathbf{x}^{i-1}_k) - f(\mathbf{x}_k')\right)}{\hat{V}}\\
    \widehat{SobT_i} &= \frac{\frac{1}{N} \sum^N_{k = 1} \left(f(\mathbf{x}^{i}_k )- f(\mathbf{x}_k')\right)^2}{2\hat{V}},
\end{split}
\end{align}
where $\mathbf{x}^*_k = (x^*_{k,1},\cdots, x^*_{k,p})$ is the $k-$th Monte-Carlo trial in the sample $\mathbf{x}^*$, $k\in \{1,n\}$ and  $\hat{V}$ is the total variance estimate that can be computed as the average of the total variances computed with each sample $\mathbf{x}^*$.

\section{Central Limit Theorem for Sobol’ indices}\label{Appendix:TCL}
We recall the theorem \ref{theo:TCL} we presented in Section \ref{section:GSA}.
\begin{theorem}
Each index $\mathcal{S}$ in the equations $\eqref{eq:defSobol’2}$ to $\eqref{eq:defSobol’Tind}$ can be written as $A/B$ and the corresponding estimate $\mathcal{S}_n$ can be written as $A_n/B_n$. For each of these indices, we have a central limit theorem:
\begin{equation}
    \sqrt{n}(\mathcal{S}_n - \mathcal{S}) \xrightarrow{D} \mathcal{N}(0, \sigma^2)
\end{equation}
with $\sigma^2$ depending on which index we study.
\end{theorem}

We propose to study the central limit theorem for the estimator of the index $Sob_i$ proposed in Appendix \ref{Appendix:estim_mara}. Note that the result is the same for other estimators of the Sobol’ indices proposed in the same section.

If we denote
\begin{equation}
  Z_n = \begin{pmatrix}
    n^{-1} \sum f(X_{i,k},X_{\sim i,k})f(X_{i,k},X^\prime_{\sim i,k}) \\
    n^{-1} \sum f(X_{i,k},X_{\sim i,k})f(X^\prime_{i,k},X^\prime_{\sim i,k})\\
    n^{-1} \sum f(X_{i,k},X_{\sim i,k})\\
    n^{-1} \sum f^2(X_{i,k},X_{\sim i,k})
  \end{pmatrix}
\end{equation}
then the estimator $\widehat{Sob_i}$ of the Sobol’ index $Sob_i$ is equal to $h(Z_n)$ where
\begin{equation*}
  h(\beta_1,\beta_2,\beta_3,\beta_4) = \frac{\beta_1 - \beta_2}{\beta_4 - \beta^2_3}.
\end{equation*}

Applying the delta-method \cite{van2000asymptotic}, we obtain the convergence of $h(Z_n)$ to $h(Z) = Sob_i$
\begin{equation}
   \sqrt{n} \left(\widehat{Sob_i} - Sob_i\right) \rightarrow \mathcal{N}(0,\nabla h(\beta)\Sigma \nabla h(\beta)^T),
 \end{equation}

for which we need to compute the gradient of $h$

\begin{equation*}
  \nabla h(\beta_1,\beta_2,\beta_3,\beta_4) = \left(\frac{1}{\beta_4 - \beta^2_3}, - \frac{1}{\beta_4 - \beta^2_3}, \frac{2(\beta_1 - \beta_2)\beta_3}{(\beta_4 - \beta^2_3)^2}, \frac{-(\beta_1 - \beta_2)}{(\beta_4 - \beta^2_3)^2}\right)^T
\end{equation*}

and the correlation matrix $\Sigma$ for the variable $Z_n$ which is

\begin{equation}
\Sigma = \begin{pmatrix}
  \sigma^2_{11} & \sigma^2_{12} & \sigma^2_{13} & \sigma^2_{14}\\
  \sigma^2_{12} & \sigma^2_{22} & 0 & 0 \\
  \sigma^2_{13} & 0 & \sigma^2_{33} & \sigma^2_{34}\\
  \sigma^2_{14} & 0 & \sigma^2_{34} & \sigma^2_{44}
\end{pmatrix}
\end{equation}
where the values $\sigma^2_{ij} = Cov(Z_{i},Z_{j})$ are given as
\begin{equation}
  \begin{split}
    \sigma^2_{11}  = & \mbox{Var}(f(X,X_{\sim i})f(X,X_{\sim i}^\prime))\\
    \sigma^2_{12}  = & \E[f^2(X,X_{\sim i})f(X,X_{\sim i}^\prime)f(X^\prime,X_{\sim i}^\prime)]\\
    \sigma^2_{13}  = & \E[f^2(X,X_{\sim i})f(X,X_{\sim i}^\prime)]\\
    \sigma^2_{14}  = & \E[f^3(X,X_{\sim i})f(X,X_{\sim i}^\prime)f(X^\prime,X_{\sim i}^\prime)] - \E[f^2(X,X_{\sim i})]\E[f(X,X_{\sim i}^\prime)f(X,X_{\sim i}^\prime)]\\
    \sigma^2_{22}  = & \mbox{Var}(f(X,X_{\sim i}))^2\\
    \sigma^2_{33}  = & \mbox{Var}(f(X,X_{\sim i}))\\
    \sigma^2_{34}  = & \E[f^3(X,X_{\sim i})]\\
    \sigma^2_{44}  = & \E[f^4(X,X_{\sim i}) - \E[f^2(X,X_{\sim i})]^2.\\
  \end{split}
\end{equation}

\section{Estimation of Cramér-von-Mises indices}\label{Appendix:estimcvm}
We propose two ways of estimating the extended Cramér-von-Mises indices that we denote by $U(Y,X_i|X_{\sim i})$ defined in \eqref{eq:extended_cramer}.

The first one is to use the fact that \begin{equation}\label{equation:cramer-von-mises-independent}
\begin{split}
    U(Y,X_i|\mathbf{Z}) & = \frac{\int \E (\mbox{Var}(\E\left[\1_{Y\leq t}|X_i,\mathbf{Z}\right]|\mathbf{Z}))d\mu(t))}{\int \mbox{Var}(\1_{Y\leq t}) d\mu(t))} \\
    & = T(Y,X_i|\mathbf{Z}) \times (1 - T(Y,\mathbf{Z})).\\
\end{split}
\end{equation}
We need to estimate $T(Y,X_i|X_{\sim i})$ and $T(Y,X_{\sim i})$. Estimates for both theses quantities are taken from \cite{azadkia2019simple}.

Consider a triple of random variables $(X,Z,Y)$ and an i.i.d sample $(X_i,Z_i,Y_i)_{1\leq i \leq n}$.  For simplicity, we still suppose the random variables to be diffuse (that is without ties). The random variable $Z$ is used for the conditioning.

For each $i$, let $N(i)$ be the index $j$ such that $Z_{j}$ is the nearest neighbor of $Z_{i}$ with respect to the Euclidean distance and let $M(i)$ be the index $j$ such that $(X_{j},Z_{j})$ is the nearest neighbor of $(X_{i},Z_{i})$. Let $R_{i}$ be the rank of $Y_{i}$, that is the number of $j$ such that $Y_{j}\leq Y_{i}$.

The correlation coefficient defined in \cite{azadkia2019simple} is defined as:
\begin{equation}
    T_{n}(Y,X|Z) = \frac{\sum_{i=1}^n \left(\min\{R_{i},R_{M(i)}\} - \min\{R_{i},R_{N(i)}\}\right)}{\sum_{i=1}^n \left(R_{i} - \min\{R_{i},R_{N(i)}\}\right)}.
\end{equation}
The authors of \cite{azadkia2019simple} prove that this estimator converges almost surely to a deterministic limit $T(Y,X|Z)$ which is equal to the quantity we defined in the first section. In order to estimate the extended Cramér-von-Mises sensitivity index $CVM^{ind}_{X}$, we propose the estimator
\begin{equation}
    U_n(Y,X_i|X_{\sim i}) = T_n(Y,X_i|X_{\sim i}) \times \left(1 -T_n(Y,X_{\sim i})\right).
\end{equation}
The convergence of the estimator $U_n(Y,X_i|X_{\sim i})$ to the quantity of interest $U(Y,X_i|X_{\sim i})$ is immediate.

We propose an alternative method for the estimation of this index. We take advantage of the estimates given in \cite{azadkia2019simple} and \cite{chatterjee2020new}. We have the two following convergences almost surely:
\begin{equation}
  Q_n(Y,X|Z)= \\ n^{-2} \sum_{j=1}^n \left(\min\{R_{j},R_{M(j)}\} - \min\{R_{j},R_{N(j)}\}\right)\\ \rightarrow \int \E (\mbox{Var}(\E\left[\1_{Y\leq t}|X,Z\right]|Z))d\mu(t))
\end{equation}
\begin{equation}
  S_n(Y) = n^{-3} \sum_{j=1}^n L_j (n-L_j) \rightarrow \int \mbox{Var}(\1_{Y\leq t}) d\mu(t))
\end{equation}
where $L_j$ is the number of $k$ such that $Y_k \geq Y_j$.

\begin{proposition}[Estimator of the extended Cramér-von-Mises indices]
  The quantity defined as $\tilde{U}_n(Y, X|Z) =Q_n(Y,X|Z)/S_n(Y)$ is a consistent estimator of $U(Y,X_i|X_{\sim i})$.
\end{proposition}
The proof is obtained directly using classical probability tools.

\section{Proofs}
\subsection{Proof of Theorem~\ref{levy}}
\begin{proof}
Indeed, we can always write
\begin{equation}\label{eq:decomposition}
    p_\mathbf{X} = p_{X_i} \times p_{X_{i+1}| X_i}\times \cdots \times p_{X_{i-1}|X_{\sim (i-1)}}.
\end{equation}
Since we are back to a product of marginals, we have a hierarchical independence. We choose the cyclical hierarchy ( $X_i$, followed by $X_{i+1}|X_i$, then  $X_{i+2}|X_i,X_{i+1}$, and so on and so forth till $X_{i-1}|X_{\sim (i-1)}$ ) as we are in fact only interested in the first and the last elements of this hierarchy ( $X_i$ and $X_{i-1}|X_{\sim (i-1)}$). We can always map univariate random variables to uniform distributions by matching the quantiles by using the cumulative distribution function -- one can view this operation as hierarchical Optimal Transport, see \cite{carlier2010knothe} -- and by doing so for each variable defined above, we have the so-called Levy-Rosenblatt transform, denoted here as RT, that is:
\begin{equation}
        (X_i, (X_{i+1}|X_i), \cdots, (X_{i-1}|X_{\sim (i-1)})) \sim p_{\mathbf{X}} \xrightarrow{RT} (U^i_1,\cdots,U^i_p) \sim \mathcal{U}^p(0,1).
\end{equation}

\end{proof}
\subsection{Proof of Examples following ~\ref{def:link}}
\begin{proof}
We will show here how each definition of fairness and GSA measure presented in Table \ref{tab:Fairness_GSA} match for binary classification with $S$ binary.

\begin{enumerate}[(i)]
  \item The definition of \textit{Statistical Parity} is given by  $|\P(f(\mathbf{X}) =1|S = 1) - \P(f(\mathbf{X}) =1|S = 0) |$. For simplicity, we consider $\mbox{Var}(f(\mathbf{X})) = 1$. If we compute the Sobol’ index of the predictor $f(\mathbf{X})$ for the protected variable $S$, we obtain:
  \begin{align*}
    Sob_S(f(\mathbf{X})) & = \mbox{Var}_S(\E_{\mathbf{X}\setminus S}[f(\mathbf{X})|S])\\
    & = \E_S \E_{\mathbf{X}\setminus S}^2[f(\mathbf{X})|S] -  \E_{\mathbf{X}}[f(\mathbf{X})|S]^2\\
    & = \P(S=1)\P(f(\mathbf{X}) =1|S = 1)^2  + \P(S=0)\P(f(\mathbf{X}) =1|S = 0)^2  - \P(f(\mathbf{X}) = 1)^2 \\
    & = \P(S = 1) \P(S = 0)\times \left[ \P(f(\mathbf{X}) =1|S = 1) - \P(f(\mathbf{X}) =1|S = 0) \right]^2 \\
    & = \P(S = 1) \P(S = 0)\times DI^2.
  \end{align*}
  We see that the quantity of interest in \textit{Statistical Parity} is the same as the Sobol’ index, up to a constant depending on the proportion in each class of the protected variable.

  \item For \textit{avoiding Disparate mistreatment}, the quantity of interest is $|\P(f(\mathbf{X}) \not = Y|S = 1) - \P(f(\mathbf{X}) \not = Y|S = 0)|$. This can be obtained by replacing $f(\mathbf{X})$ by $\1_{f(\mathbf{X}) \not = Y}$ in the quantity of interest for \textit{Statistical Parity}.
  Therefore, by the same computation as previously, we can link \textit{avoiding Disparate mistreatment} to the Sobol’ index of the error of the predictor $\1_{f(\mathbf{X}) \not = Y}$ for the protected variable $S$.

  \item For \textit{Equality of Odds}, we are interest in the difference $|\P(f(\mathbf{X}) |Y = i, S = 1) - \P(f(\mathbf{X})| Y = i, S = 0)|$ for $i=0,1$. Each of this difference can be expressed as seen before as $\mbox{Var}_S(E_X[f(\mathbf{X})|Y = i, S])$. Since we want this quantity to be equal to zero for each $i$, we can compute \textit{Equality of Odds} with $\E_{Y} \mbox{Var}_S(E_X[f(\mathbf{X})|Y, S])$, which is the extended Cramèr-von-Mises index of the predictor for the protected variable $S$.

  \item For \textit{avoiding Disparate Treatment}, the quantity of interest is very similar to \textit{Statistical Parity} since we are interested in proving $f(\mathbf{X})|\mathbf{X}\setminus S \ind S$. By similar computations as before, this fairness boils back to looking at $\E_{\mathbf{X}\setminus S} \mbox{Var} \E_{\mathbf{X}\setminus S} [f(\mathbf{X})|\mathbf{X}]$. This can be simplified into $\E_{\mathbf{X}\setminus S} \mbox{Var} [f(\mathbf{X})|\mathbf{X}\setminus S]$, which is the Total Sobol’ index of the predictor for the protected variable $S$.

\end{enumerate}

\end{proof}

\subsection{Proof of Proposition~\ref{prop:causal}}
\begin{proof}
The proof is a direct consequence of the Hoeffding decomposition of the function $Y = \psi(X,S)$. By factorizing $\P_Y$ as $\P_{Y|X,S}\P_{X|S}\P_{S}$, we can write \begin{equation*}
  Y = \psi_{X}(X(S)) + \psi_{S}(S) + \psi_{S,X}(S)\times \psi_{X,S}(X(S))
\end{equation*}
If $SobT^{ind}_S = 0$ then $\mbox{Var}(\psi_{S}(S) + \psi_{S,X}(S)\times \psi_{X,S}(X(S))) = 0$. By orthogonality in the Hoeffding decomposition, $\mbox{Var}(\psi_{S}(S)) = \mbox{Var}(\psi_{S,X}(S)\times \psi_{X,S}(X(S))) = 0$, which lead to $\psi_{S}(S) = \psi_{S,X}(S)\times \psi_{X,S}(X(S)) = 0$. It holds that $Y = \psi_{X}(X(S))$.

For the second part of the proposition, we apply the same reasoning by factorizing $\P_Y$ as $\P_{Y|X,S}\P_{S|X}\P_{X}$.
We can write \begin{equation*}
  Y = \psi^\prime_{S}(S(X)) + \psi^\prime_{X}(X) + \psi^\prime_{S,X}(X)\times \psi^\prime_{X,S}(S(X))
\end{equation*}
If $SobT_S = 0$ then $\mbox{Var}(\psi^\prime_{S}(S(X)) + \psi^\prime_{S,X}(X)\times \psi^\prime_{X,S}(S(X))) = 0$. By orthogonality in the Hoeffding decomposition, $\mbox{Var}(\psi^\prime_{S}(S(X))) = \mbox{Var}(\psi^\prime_{S,X}(X)\times \psi^\prime_{X,S}(S(X)) = 0$, which lead to $\psi^\prime_{S}(S(X)) = \psi^\prime_{S,X}(X)\times \psi^\prime_{X,S}(S(X)) = 0$. It holds that $Y = \psi^\prime_{X}(X)$.

\end{proof}

\subsection{Proof of Proposition~\ref{prop:intersect1} and Proposition~\ref{prop:intersect2}}

\begin{proof}
Without loss of generality, we can consider only two sensitive features $S_1$ and $S_2$. Because of the various bounds on Sobol’ indices explained in previous Section, we know that $SobT_{S_1, S_2} \leq SobT_{S_1}$. $SobT_{S_1}$ is the GSA measure associated with \textit{Avoiding Disparate Treatment}. This means that to be fair in the sense of \textit{Avoiding Disparate Treatment} implies the nullity of $SobT_{S_1}$ and therefore the nullity of $SobT_{S_1, S_2}$. The second result is a direct consequence of the absence of bounds between $Sob_{S_1}$ and Sobol’ indices for $(S_1, S_2)$ and an example has been given in the previous toy-case in introduction of the Subsection. We can find cases where $Sob_{S_1}$ is arbitrary high and $Sob_{S_1,S_2}$ is null, such as $f(X) = S_1$; and cases where $Sob_{S_1}$ is null and $Sob_{S_1,S_2}$ is arbitrary high, such as $f(X) = S_1 \times S_2$.
\end{proof}

\end{document}